\documentclass[12pt]{article}
\usepackage{fullpage,graphicx,amsmath,amsfonts}
\usepackage[small,bf]{caption}
\usepackage{subcaption}
\usepackage{enumitem}
\usepackage{authblk}
\usepackage{url}
\usepackage{booktabs}
\usepackage{hyperref}
\usepackage{multirow}

\hypersetup{
    colorlinks=true,     
    linkcolor=blue,      
    citecolor=red,       
    urlcolor=magenta     
}

\newcommand{\reals}{{\mbox{\bf R}}}

\newcommand{\eg}{{\it e.g.}}
\newcommand{\ie}{{\it i.e.}}

\newcommand{\BEAS}{\begin{eqnarray*}}
\newcommand{\EEAS}{\end{eqnarray*}}
\newcommand{\BEA}{\begin{eqnarray}}
\newcommand{\EEA}{\end{eqnarray}}
\newcommand{\BEQ}{\begin{equation}}
\newcommand{\EEQ}{\end{equation}}
\newcommand{\BIT}{\begin{itemize}}
\newcommand{\EIT}{\end{itemize}}
\newcommand{\BNUM}{\begin{enumerate}}
\newcommand{\ENUM}{\end{enumerate}}

\newcommand{\BEQN}{\begin{equation*}}
\newcommand{\EEQN}{\end{equation*}}

\newcounter{algorithmctr}
\renewcommand{\thealgorithmctr}{\arabic{algorithmctr}}
   {\mbox{}\\*[\parskip]\begin{minipage}{\linewidth}%
       \refstepcounter{algorithmctr}\begin{list}{}{%
       \setlength{\rightmargin}{0\linewidth}%
       \setlength{\leftmargin}{.05\linewidth}}%
       \rmfamily\small
       \item[]{\setlength{\parskip}{0ex}\hrulefill\par%
        \nopagebreak{\bfseries\textsf{Algorithm \thealgorithmctr~}}}}%
   {{\setlength{\parskip}{-1ex}\nopagebreak\par\hrulefill\\*[2ex]\par}%
   \end{list}\end{minipage}}

\title{Presolving for GPU-Accelerated First-Order LP Solvers}
\author[]{Daniel Cederberg \qquad \qquad Stephen Boyd \\ Stanford University}

\begin{document}
\maketitle

\begin{abstract}
Recent research has focused on developing GPU-accelerated first-order solvers
for linear programming (LP). This line of work, however, has largely overlooked
the role of presolving, and thus prior results do not fully reflect the speedups
achievable through GPU acceleration in a realistic end-to-end solver pipeline.
At the same time, LP presolving has traditionally been developed for CPU-based
solvers, where presolve time rarely dominates the total runtime and the emphasis
has been on maximizing the reduction in problem size, even at the expense of
costly presolve rules. Given the high performance of modern GPU-accelerated
solvers and the inherently sequential nature of presolving, it is unclear
whether this traditional approach to presolving remains appropriate.

In this paper we revisit LP presolving from the perspective of GPU-accelerated
first-order LP solvers. We identify a set of relatively simple presolve rules and
show that a carefully engineered collection of these captures most of the
reduction achieved by Gurobi's commercial state-of-the-art presolver, at a
fraction of the cost. Moreover, we demonstrate that such lightweight presolving
can yield substantial end-to-end speedups for the GPU-accelerated solver
cuPDLPx, despite presolving sometimes constituting a significant fraction of the
total runtime. These results suggest that lightweight presolving may remain
beneficial as GPU performance continues to scale, while the sequential nature of
presolving presumably does not.

We accompany this paper with an open-source C implementation of an LP presolver,
called PSLP (Presolver for Linear Programs). PSLP is battle-tested and has been
adopted by the community, with integrations in cuPDLPx, cuOpt (NVIDIA's
optimization library), and HPR-LP.
\end{abstract}
      
\newpage
\tableofcontents
\newpage
\section{Introduction}
Linear programming (LP) is arguably the most widely used problem class in
mathematical optimization, with applications spanning a broad range of fields.
Decades of research and substantial engineering effort have produced highly
optimized implementations of the simplex method and its many variants
\cite{Dantzig1951, Maros2002}, as well as interior-point methods
\cite{Roos2005}. As a result, LP solving has become a mature and remarkably
efficient technology --- capable of handling instances of a scale and complexity
that George Dantzig himself could hardly have imagined, often in a matter of
seconds.

Commercial LP solvers consist of three major components: (1) a \emph{presolver}
that reduces the size of the problem and later reconstructs a solution to the
original problem \cite{Andersen1995a}; (2) a \emph{core step} that runs one or
multiple algorithms concurrently to solve the reduced problem; and (3) a
\emph{crossover mechanism} (used only for non-simplex algorithms) that takes the
approximate solution produced by the core step and converts it into a highly
accurate basic feasible solution \cite{Bixby1994}. While the first and third
components are highly sequential in nature (and can take a significant portion
of the total solve time), the core step traditionally implements an
interior-point method that benefits from CPU multithreading for the
computational bottleneck of factorizing and solving large sparse linear systems.
However, despite being compute-bound, this bottleneck appears to have an
inherent limit to the amount of parallelism that can be effectively exploited
(see, for example, \cite{Rothberg2020, Panos2022}).

To fully exploit the computational power of modern hardware accelerators in the
context of linear programming, recent research has focused on speeding up the
core step by developing first-order methods that rely heavily on
matrix-vector multiplications \cite{Applegate2021, Applegate2025, Chen2025,
Lu2023, Lu2025}. While this paradigm shift opens up exciting opportunities for
solving extremely large LPs, it also raises important questions about the
role of the other traditional components of LP solvers.  In this paper we
focus on presolving, a critical but often overlooked component, and its
impact on GPU-accelerated first-order solvers.

Presolving is traditionally performed on CPUs because it involves a large number
of sequential, logic-driven, and structurally irregular problem transformations.
The application of one problem transformation often enables or invalidates
others, creating a chain of data-dependent transformations that are inherently
sequential. Many of these transformations require heavy branching, dynamic
sparse-matrix modifications, and unpredictable memory access patterns. It is
therefore currently believed that CPUs, with their low-latency caches and
sophisticated control-flow capabilities, are far better suited for the task than
GPUs. Although certain parts of a presolver admit some degree of parallelism,
the engineering effort required to map these onto GPU architectures is rarely
justified, and we speculate that the potential speedup would at best be modest,
in particular if CPU multithreading is used for the parallelizable parts.
However, GPU acceleration has a different consequence for presolving: the
relative cost of presolving compared to the GPU-accelerated core step has
increased significantly, putting very high demands on the implementation and
speed of the presolver. For every extra second spent in presolving, there are
thousands of idle GPU cores waiting to get to work, so it is unclear whether LP
presolving is still worthwhile in this new context. 

The idea of reducing the size of an LP before passing it to a solver has a long
history, dating back at least to the 1970s. (We review related work in \S
\ref{sec:related-work}.) Over the past decades, a wide range of presolve
techniques have been proposed, from simple ones such as eliminating
constraints with just one variable \cite{Brearley1975} to more elaborate
procedures such as exploiting symmetry between variables \cite{Achterberg2019}.
This body of work has largely been developed for CPU-based solvers, where
presolving and the core algorithm run on the same hardware, and where presolve
time rarely dominates the overall runtime. In this traditional setting, the
emphasis has been on maximizing the reduction in problem size, even at the
cost of more costly presolve techniques.

In this paper we revisit LP presolving from a different perspective, motivated
by GPU-accelerated first-order methods. In this setting, presolving is
inherently sequential, runs on different hardware than the core algorithm, and
can easily become a bottleneck. This raises a fundamental question: can
meaningful performance gains be obtained using only a carefully engineered
collection of simple and fast presolve techniques, and how much problem
reduction is sacrificed by forgoing the more complex and costly rules found in
traditional presolvers? By focusing on rules that are fast and lightweight, we
aim to identify those that have the potential to remain effective as GPU
performance continues to scale, while the sequential nature of presolving
presumably does not.

Our experiments show, perhaps surprisingly, that a carefully engineered
collection of relatively simple presolve techniques can achieve reductions in
problem size comparable to those obtained by the commercial state-of-the-art
presolver implemented in Gurobi \cite{Achterberg2019}. Specifically, we find
that our open-source C implementation of an LP presolver, called PSLP (Presolver
for Linear Programs), on average captures about 90\% of the reduction achieved
by Gurobi's presolver, while being over 6 and 11 times faster, respectively, on
two standard LP datasets. We also investigate the impact of PSLP on the total
solve time of cuPDLPx \cite{Lu2025-cuPDLPx}, a GPU-accelerated first-order LP
solver based on the primal-dual hybrid gradient (PDHG) method \cite{Zhu2008,
Chambolle2011}. Our main observation is that PSLP has a strong positive impact
on the total solve time of cuPDLPx across a wide range of LP instances. For
example, on the large problems with more than $10^7$ nonzero entries from
Mittelmann's LP collection \cite{mittelmannLP} and the MIPLIB 2017 root-node LP
relaxations \cite{Gleixner2021}, enabling PSLP reduces the shifted geometric
mean of the total solve time by more than a factor of 2.5 and 9, respectively.
Since PSLP was released, it has been integrated into several solvers by their
developers, including cuPDLPx, cuOpt (NVIDIA's GPU-accelerated optimization
library \cite{CuOpt2025}), and HPR-LP \cite{Chen2025}.

The remainder of this paper begins with an overview of related work on
presolving. Next, we review the main ideas behind LP presolving in \S
\ref{sec:background}. Our presentation emphasizes the abstractions of
\emph{reductions}, \emph{primal exploration}, and \emph{dual exploration}, which
we believe are useful for understanding and implementing presolvers, although
the classic literature never presents presolving in this way. In \S
\ref{sec:design-and-implementation}, we discuss some of the design choices
behind our open-source implementation of a fast and lightweight LP presolver. In
\S \ref{sec:experiments}, we assess how much of the total reduction achieved by
the commercial presolver implemented in Gurobi can be captured by the carefully
engineered collection of relatively simple presolve techniques implemented in
PSLP. We also evaluate the impact of presolving on the total solve time of
cuPDLPx. Finally, we discuss our findings in \S \ref{sec:discussion}.

\subsection{Related work}
\label{sec:related-work}
The idea of reducing the size of an optimization problem before passing it to a
solver has a long history. Early references describing many of the transformations
implemented in modern presolvers are \cite{Brearley1975, Williams1983,
Bradley1983, Tomlin1986}. These references focus primarily on linear programming
where the simplex method is used as the algorithm for the core step. The
development of interior-point methods in the 90's sparked renewed interest in
presolving, as evidenced by a string of papers investigating the importance of
presolving for interior-point methods for LPs \cite{Andersen1995a,
Andersen1995b, Gondzio97}. An efficient LP presolver forms the foundation of
presolvers for more general problem classes, including quadratic programs
\cite{Meszaros2003, Gould2004}, mixed-integer linear programming (MIP)
\cite{Gamrath2015, Achterberg2019}, and nonlinear programming \cite{Zhang2025}.

Despite the well-established role of presolving in commercial optimization
software, there are surprisingly few open-source presolvers of any kind. One
notable exception is PaPILO \cite{Gleixner2023}, a solver-independent presolver
for MIP developed by Leona Gottwald with significant contributions from Ivet
Galabova, Alexander Hoen, Ambros Gleixner, and several others (see the PaPILO
GitHub page for the full list of contributors). The scarcity of open-source
presolvers is likely due to the substantial engineering effort required to
implement them: presolvers involve numerous subtleties and corner cases that
make them far more complex to implement than their high-level descriptions
suggest. Moreover, a slow and poorly implemented presolver can do more 
harm than good and will not be used. 
In addition, there is relatively little academic
incentive to develop presolvers, since presolvers are often viewed as
engineering components rather than novel research contributions. This lack of
incentive is reflected in the number of modern open-source solvers developed in
academia that ship without a presolver, including Clarabel \cite{Goulart2024},
SCS \cite{ODonoghue2016}, OSQP \cite{Stellato2020}, ECOS \cite{Domahidi2013},
cuPDLP \cite{Lu2023, Lu2025}, cuPDLPx \cite{Lu2025-cuPDLPx} (that now ships with
PSLP), and several others. Notable exceptions are HiGHS
\cite{Huangfu2018, Galabova2023}, CBC \cite{Forrest2005}, and Google's linear
optimization package \cite{ORTOOLS2025}, which all include solver-dependent
presolve functionality for LPs.

Finally, we mention that while we view presolving as a task that runs on a CPU,
there has been work on accelerating the so-called \emph{primal propagation} or
\emph{bound tightening} step of a presolver using a GPU \cite{Sofranac2022}.
While the results in that paper are promising given the setup and performance
metrics they use, it is unclear how their results would translate to a more
realistic setting for LP presolving for the following two reasons. First, they
run up to 100 rounds of primal propagation, which is far more than what is
typically done in an LP presolver in practice. (For example, Gurobi's presolver
\cite[\S3.2]{Achterberg2019} and PSLP both propagate each constraint only
\emph{once} per round, and typically run only for a few rounds. Moreover, PSLP
keeps track of which constraints need to be re-examined after each round, so in
every round only a subset of the constraints are actually propagated.) Second,
primal propagation is only one of many components of a presolver. For example,
in our experiments on Mittelmann's LP collection \cite{mittelmannLP}, primal
propagation on average accounts for 12.0\% of the total time spent in PSLP.

\section{Background}
\label{sec:background}
We consider an LP of the form 
\begin{equation}\label{eq:lp}
\begin{array}{ll}
\mbox{minimize} & c^T x \\
\mbox{subject to} & \underline{b} \leq Ax \leq \overline{b} \\
& \underline{x} \leq x \leq \overline{x}, \\
\end{array}
\end{equation}
with variable $x \in \reals^n$. The problem data are the cost vector $c \in
\reals^n$, the left and right hand sides $\underline{b} \in \reals^m$ and
$\overline{b} \in \reals^m$, the constraint matrix $A \in \reals^{m \times n}$,
and the variable bounds $\underline{x} \in \reals^n$ and $\overline{x} \in
\reals^n$. Infinite values are allowed in both the variable bounds and the left and
right hand sides, allowing for unbounded variables and one-sided inequality
constraints. Equality constraints can be expressed by setting the corresponding
components of $\underline{b}$ and $\overline{b}$ equal.

A solution of \eqref{eq:lp} consists of a primal solution $x^\star$
and a dual solution $(y^\star, z^\star)$ satisfying the optimality conditions
\begin{equation} \label{eq:KKT}
\begin{split}
& \underline{b} \leq Ax^\star  \leq \overline{b}, \qquad
\underline{x} \leq x^\star \leq \overline{x} \\
& c - A^T y^\star - z^\star  =  0 \\
& z_k^\star > 0 \implies x_k^\star  = \underline{x}_k \\
& z_k^\star < 0 \implies x_k^\star  = \overline{x}_k \\
& y_i^\star > 0 \implies (Ax^\star)_i  = \ \underline{b}_i \\
& y_i^\star < 0 \implies (Ax^\star)_i  = \ \overline{b}_i.
\end{split}
\end{equation}
A few subtleties are worth pointing out. First, for an unbounded variable $x_k$
(with $\underline{x}_k = -\infty$ and $\overline{x}_k = +\infty$), the
corresponding dual variable $z_k^\star$ must be zero. Second, if
$\underline{b}_i = -\infty$, then $y_i^\star \leq 0$, and if $\overline{b}_i =
+\infty$, then $y_i^\star \geq 0$. Third, if there does not exist any
primal-dual pair $(x^\star, y^\star, z^\star)$ satisfying these conditions, then
the problem is either infeasible or unbounded.

In the following subsections, we summarize some of the most important ideas of presolving an
LP of the form \eqref{eq:lp}. First we introduce the notion of \emph{reductions}, each
of which consists of a \emph{presolve transformation} and a \emph{postsolve
transformation}. Next, we describe \emph{primal} and \emph{dual
exploration}, which analyze the primal and dual constraints to identify
opportunities for invoking reductions.

\subsection{Reductions}
The idea of a presolver is to apply a sequence of \emph{reductions} to the
original problem to obtain a smaller, equivalent problem that can be solved more
efficiently. The notion of reductions is standard in both the theory and
practice of computer science, where it underlies, for example, program rewriting
systems in compilers \cite{aho2006} and transformation-based modeling frameworks for
optimization such as CVXPY \cite{agrawal2018rewriting}.

In the context of LP presolving, each reduction consists of a \emph{presolve
transformation} and a \emph{postsolve transformation}. The presolve
transformation takes as input one problem and outputs an equivalent problem that
is, in some sense, smaller or simpler than the input problem. The postsolve
transformation takes as input a solution to the reduced problem and maps it to a
solution of the original problem. Each reduction is very simple and cheap to
apply, and for large LPs several hundreds of thousands of reductions may be
applied in sequence. 

A first example of a reduction is \emph{fixing a variable}. Suppose we for some
reason fix a variable $x_k$ to some value $\hat{x}_k$ (we discuss how to
identify such variable fixings below), and assume that this is the only
reduction we apply. The \emph{presolve} transformation for this reduction
consists of removing the variable $x_k$ from the original problem by
substituting its fixed value into the remaining parts of the problem. The
variable $x_k$ appears in the original problem but not in the reduced problem,
so the \emph{postsolve} transformation must recover an optimal primal variable
$(x_k^\star)_{\text{original}}$ and an optimal dual variable
$(z_k^\star)_{\text{original}}$ of the original problem from a solution of the
reduced problem. The \emph{primal} postsolve transformation simply sets
$(x_k^\star)_{\text{original}} = \hat{x}_k$. The \emph{dual} postsolve
transformation sets $(y^\star)_{\text{original}} = (y^\star)_{\text{reduced}}$ and
$(z_k^\star)_{\text{original}} = c_k - a_k^T
(y^\star)_{\text{reduced}}$, where $a_k$ is the $k$th column of $A$ and
$(y^\star)_{\text{reduced}}$ is the optimal dual variable of the reduced
problem. (This choice of $(z_k^\star)_{\text{original}}$ follows from 
the dual feasibility condition, \ie, the second line of \eqref{eq:KKT}.)

In addition to the reduction of fixing a variable, there are four 
main reductions:
\begin{itemize}
\item removing a constraint,
\item adding a multiple of an equality constraint to another constraint,
\item substituting a variable appearing in only one equality constraint, and
\item changing variable bounds.
\end{itemize}
While classic papers on presolving, such as \cite{Brearley1975, Andersen1995a,
Achterberg2019}, list several more complicated problem transformations as
reductions, we use a different abstraction where we break down complex
transformations into sequences of simple reductions. In our abstraction, a
reduction is the \emph{simplest} atomic rewriting unit of a problem. Most
problem transformations described in the literature on LP presolve can be broken
down into sequences of the 5 reductions mentioned above. This abstraction does
not only make it easier to describe a presolver, but it also simplifies the
implementation of one of the most complicated parts of the presolver: recovering
optimal dual variables of the original problem from an optimal primal-dual
solution of the reduced problem.


To find valid reductions (\ie, problem transformations that reduce the size of
the problem while provably maintaining an equivalence to the original problem),
we use \emph{primal exploration} and \emph{dual exploration}, as described next.

\subsection{Primal exploration}
Primal exploration analyzes the \emph{primal} constraints to identify
opportunities for reductions. The logical rules applied during this process
guarantees that the feasible set of the reduced problem is equivalent to that of
the original formulation. Roughly speaking, primal exploration identifies
reductions that preserve the feasible set while producing a more compact
representation with less redundancy.

A simple example of primal exploration is the \emph{redundant constraint
explorer}, which uses variable bounds to identify redundant constraints. For
example, consider the bounds $x_1 \geq 1$ and $x_2 \geq 2$ together with the
constraint $x_1 + 2x_2 \geq 5$. Since the variable bounds imply that the
constraint is satisfied, the constraint is classified as redundant. Once such a
constraint is identified, the redundant constraint explorer invokes the
corresponding reduction to remove it from the problem.

Another example is the \emph{doubleton row explorer}, which finds equality
constraints with only two variables and substitutes one variable in terms of the
other. For example, consider the constraint $x_1 + x_2 = 1$ together with the
bounds $x_1, x_2 \geq 0$. First, we transfer the bound $x_1 \geq 0$ to an upper
bound on $x_2$, giving $x_2 \leq 1$ and making $x_1$ unbounded. (This is done by
invoking the reduction that changes variable bounds.) Second, we substitute $x_1
= 1 - x_2$  in the remaining parts of the problem. (This is done by first
repeatedly invoking the reduction that adds a multiple of an equality constraint
to another constraint, followed by the reduction that substitutes a variable
appearing in only one equality constraint.) Third, we remove the constraint $x_1
+ x_2 = 1$. (This is done by invoking the reduction that deletes a constraint.)
Once the doubleton row explorer identifies an equality constraint with only two
variables, it invokes these reductions in sequence. The reduced problem has one
less variable and one less constraint, and the reduced feasible set is
equivalent to the feasible set of the original problem in the sense that if $x_2
= a$ for some $a \in \reals$ is feasible for the reduced problem, then $(x_1,
x_2) = (1 - a, a)$ is feasible in the original problem.

We list more examples of primal exploration in \S 3.1.

\subsection{Dual exploration}
Dual exploration analyzes the \emph{dual} constraints to identify opportunities
for reductions based on complementary slackness and dual variable bounds. Unlike
primal exploration, dual exploration may identify reductions that do not
preserve the feasible set of the original problem. However, the logical rules
applied during dual exploration guarantee that any optimal solution of the
reduced problem can still be mapped to an optimal solution of the original
problem.

Dual exploration is further categorized into \emph{strong} and \emph{weak}
exploration. The semantics behind weak dual exploration imply that any reduction it
identifies preserves the entire optimal solution set: \emph{every} optimal
solution of the original problem can, in principle, be recovered from some optimal
solution of the reduced problem. In contrast, the semantics of strong dual
exploration are looser: it may yield more aggressive reductions, but these only
guarantee that \emph{at least one} optimal solution of the original problem can be
recovered from an optimal solution of the reduced problem.

A simple example of dual exploration is the \emph{variable lock explorer}. An
\emph{uplock} of a variable is a constraint that prevents the variable from
increasing to infinity. If a variable $x_k$ has objective coefficient $c_k < 0$
and no uplock, then the dual constraints imply that $z_k < 0$ for any dual
feasible point. By complementary slackness, this means that $x_k$ will be equal
to its upper bound $\overline{x}_k$ at any optimal solution, so the variable
lock explorer can invoke the reduction that fixes $x_k$ to $\overline{x}_k$. (If
$\overline{x}_k = +\infty$, then the problem is unbounded.) This is an example
of weak dual exploration, since \emph{any} optimal solution of the original
problem must have $x_k = \overline{x}_k$. If instead $c_k = 0$ and $x_k$ has no
uplock, the variable lock explorer can still invoke the reduction of fixing
$x_k$ to $\overline{x}_k$, but not every primal optimal solution $x^\star$ of
the original problem is guaranteed to have $x_k^\star = \overline{x}_k$, so this
is an example of strong dual exploration.

The variable lock explorer is a special case of a more general dual exploration
phase that propagates the dual constraints to strengthen dual variable bounds
and identify reductions based on complementary slackness (see, \eg, \cite[\S
3.4]{Andersen1995a}). Other examples of dual exploration include searching for
\emph{column singletons in inequality constraints} \cite[\S 1.2]{Brearley1975}
and \emph{dominated columns} \cite[\S 4]{Gamrath2015}. Interestingly, many dual
exploration rules can be interpreted as primal explorers applied to the dual
problem. For example, both the variable lock explorer and the dominated columns
explorer are just special primal redundant constraint explorers applied to the dual
problem. 

\section{Design and implementation}
\label{sec:design-and-implementation}
While the main concepts behind LP presolving are straightforward and simple to
describe, implementing them in practice involves numerous design choices. In
this section, we discuss some of the most important decisions and describe 
our implementation.

\subsection{Explorers}
The most important design choice is which explorers to implement. Each explorer
applies a set of logical rules to identify opportunities for reductions. Adding
more explorers typically yields smaller reduced problems, but at the cost of
increased presolve time. Table \ref{tab:explorers} describes and groups commonly
used explorers into three categories based on their computational cost. 
\begin{itemize}
\item
\emph{Fast} explorers require minimal computation to identify reductions, as
they operate primarily by scanning \emph{internal statistics} maintained by the
presolver (we describe what internal statistics are in greater detail in \S 3.2).
\item
\emph{Medium} explorers require additional computation and cannot be executed
using internal statistics alone. For example, the explorer that detects parallel
rows belongs to this category because it requires a two-level hashing algorithm
together with pairwise comparison of rows to identify constraints that are
parallel \cite[\S 5.2]{Achterberg2019}. 
\item
\emph{Slow} explorers require substantially more computation. For example, one
explorer in this category detects linear dependence among equality constraints
and requires rank-revealing matrix factorizations \cite{Andersen1995b,
Miranian2003}. A related explorer identifies opportunities for reducing the
number of nonzero entries in the constraint matrix by adding a multiple of an
equality constraint to another constraint to make the latter sparser
\cite{Gondzio97}, similar to Gaussian elimination.
\end{itemize}
The categorization in table \ref{tab:explorers} is based on our experience with
presolving and is therefore somewhat subjective. 

Almost all of the explorers listed in table \ref{tab:explorers} were proposed
over several decades by researchers with either CPU-based simplex or
interior-point solvers in mind, so it is natural to ask whether all of them are
worthwhile for GPU-accelerated first-order solvers. For example, first-order
solvers open up new opportunities for solving extremely large LPs where
traditional simplex and interior-point solvers may be too slow or run out of
memory when factorizing large sparse matrices. Since the linear dependence
explorer requires matrix factorization itself, it seems counterintuitive to
include it in the presolve for GPU-accelerated first-order solvers, as one of
the main motivations for using such solvers is to \emph{avoid} matrix
factorization. Moreover, our experiments in \S \ref{sec:experiments} show that
most of the problem reduction (on average about 90\%) is detected by fast and
medium explorers. It is therefore unclear whether the additional reductions
detected by slow explorers generally justify the extra presolve time for
problems where first-order methods may be the only viable option.

\begin{table}
\centering
{\footnotesize
\begin{tabular}{p{4.5cm} p{10cm}}
\toprule
\textbf{Explorer name} & \textbf{Description} \\
\midrule
\multicolumn{2}{c}{\textit{Fast explorers}} \\
\hline \hline
\addlinespace
Singleton rows \hfill (P) & Searches for rows with a single nonzero entry, which
 can be used to fix the value or update the bounds of a variable
 \cite[\S1.2]{Brearley1975}. \\
\hline
Doubleton rows \hfill (P) & Searches for equality constraints with two nonzero
entries, allowing one variable to be eliminated via substitution \cite{Bradley1983}. \\
\hline
Redundant constraints \hfill (P) & Searches for constraints that are implied by
variable bounds and can thus be deleted \cite[\S 1.2]{Brearley1975}. \\
\hline
\parbox[t]{4.5cm}{%
Column singleton \hfill (P)\\[-0.3ex]
{\footnotesize (in equality constraint)} } & Searches for equality constraints
with a variable appearing in only that constraint, possibly allowing the variable to be substituted and the
constraint removed \cite[\S3.2]{Andersen1995a}. \\
\hline
Variable locks \hfill (D) & Searches for variables whose locks allow them to be
fixed \cite[\S 4.4]{Achterberg2019}. \\
\hline
\parbox[t]{4.5cm}{%
Column singleton \hfill (D)\\[-0.3ex]
{\footnotesize (in inequality constraint)} } & Searches for inequality
constraints with a singleton column, possibly allowing the variable to be fixed
or the  inequality constraint converted to an equality constraint \cite[\S 1.2]{Brearley1975}.\\
\midrule

\multicolumn{2}{c}{\textit{Medium explorers}} \\
\hline \hline
\addlinespace
Parallel rows \hfill (P) & Searches for constraints that are parallel, allowing
many of them to be deleted \cite{Tomlin1986}, \cite[\S 5.2]{Achterberg2019}. \\
\hline
Parallel columns \hfill (P) & Searches for columns that are parallel, allowing
many of them to be aggregated \cite[\S 3.6]{Andersen1995a}, \cite[\S
6.3]{Achterberg2019}.  \\
\hline
Primal propagation \hfill (P) & Searches for variable bounds implied by the
constraints \cite[\S1.2]{Brearley1975}, \cite[\S 3.2]{Achterberg2019}. \\
\hline
Dual propagation \hfill (D) &  Searches for dual variable bounds implied by the
dual constraints \cite[\S 3.4]{Andersen1995a}, \cite[\S 7.5]{Achterberg2019}. \\
\midrule

\multicolumn{2}{c}{\textit{Slow explorers}} \\
\hline \hline
\addlinespace
Linear dependence  \hfill (P) & Searches for linearly dependent equality
constraints \cite{Andersen1995b, Miranian2003}. \\
\hline
Variable substitution \hfill (P) & Searches for variables that can be
eliminated via substitution without introducing excessive fill-in \cite[\S 2.2]{Meszaros2003},
\cite[\S 4.5]{Achterberg2019}. \\
\hline
Variable symmetry \hfill (P) & Searches for symmetries among variables, allowing
variables to be aggregated \cite{grohe2014dimension}, \cite[\S
7.1]{Achterberg2019}. \\
\hline
Primal sparsification \hfill (P) & Searches for opportunities to combine
constraints to make them sparser \cite{Gondzio97}, \cite[\S
5.3]{Achterberg2019}. \\
\hline
Dual sparsification \hfill (D) & Searches for opportunities to combine dual
constraints to make them sparser \cite{Gemander2020}. \\
\hline
Dominated columns \hfill (D) & Searches for pairs of variables where one is
\emph{dominated}, \ie, worse than the other with respect to both the constraints and
the objective, and can thus possibly be fixed \cite[\S 4]{Gamrath2015}, \cite[\S
6.4]{Achterberg2019}. \\
\bottomrule
\end{tabular}
} \caption{Explorers grouped by computational cost. The letters (P) and (D)
indicate whether the explorer is primarily primal or dual.}
\label{tab:explorers}
\end{table}

\subsection{Data structures}
While most exploration rules are simple to describe, their efficient
implementation is often quite involved. To make the exploration phase efficient,
a presolver keeps track of several internal statistics, which typically include
the number of nonzeros in each row and column of $A$, the variable locks (\ie, the
number of constraints that prevent a variable from increasing or decreasing to
$\pm \infty$), and the smallest and largest value a constraint can take given
the current variable bounds (these are often referred to as the \emph{minimum}
and \emph{maximum activity} of a constraint \cite[\S 2]{Achterberg2019}). The
internal statistics must be maintained, so they are updated each time a
reduction is applied. Furthermore, presolving is arranged in \emph{rounds}. To
avoid useless work from re-examining all constraints in every round to identify
opportunities for reductions, one option is to keep track of the constraints
modified in the previous round and apply operations such as primal propagation
only to those constraints.

Some of the explorers require fast access to rows while others require fast
access to columns. One option is to store the constraint matrix $A$ in
compressed sparse row (CSR) format to enable fast and cache-friendly access to
rows, and supplement the structure of $A$ by a linked lists superstructure
describing its columns \cite{Gould2004}. Another option, adopted by PaPILO
\cite{Gleixner2023} and our presolver PSLP, is to store $A$ in both CSR and
compressed sparse column (CSC) formats. Storing $A$ in CSC format enables fast
and cache-friendly access to columns, but care must be taken to ensure that the
two representations of $A$ remain consistent once reductions are applied.
Furthermore, some reductions require modifying the sparsity pattern of $A$ and
may lead to a larger number of nonzero entries in certain rows or columns, so it
is useful to append each row in the CSR representation and each column in the
CSC representation with a small amount of extra space to allow for new nonzero
entries. If the extra space is exhausted for a specific row or column, nearby
rows or columns may have more extra space, which can be used after shifting the
data for nearby rows or columns.

\subsection{Our implementation}
Our implementation, which is based on the reduction-based abstraction described
in \S \ref{sec:background}, is called PSLP (Presolver for Linear Programs) and
is available at
\begin{center}
\url{https://github.com/dance858/PSLP}.
\end{center}
PSLP is implemented in C and runs (mostly) single-threaded on the CPU, with
certain parts parallelized on two or four POSIX threads. To avoid becoming a
bottleneck when used together with a GPU-accelerated LP solver, PSLP is designed
to be fast and lightweight. Consequently, it implements only the fast and medium
explorers listed in table \ref{tab:explorers}, and omits the slow explorers. For
certain problem instances, the additional reductions detected by slow explorers
could be large enough to justify their use. However, PSLP does not currently
implement slow explorers due to (1) the substantial additional engineering
effort required to implement them efficiently, and (2) their marginal benefit in
terms of problem-size reduction appears limited. Specifically, PSLP already
attains, on average, about a 90\% reduction in problem size relative to Gurobi's
presolver that implements the slow explorers \cite{Achterberg2019}, while being
significantly faster (see table \ref{tab:presolver_comparison} in \S
\ref{sec:experiments}).

PSLP is solver independent, meaning that it can be used as a standalone
presolver module for any LP solver. Its interface is minimal and just requires
the user to provide pointers to the problem data, so solvers that lack a
built-in presolver can easily integrate it as an external component with just a
few lines of code. PSLP takes in an LP of the form \eqref{eq:lp}, and emits a
reduced LP on the same form. Given an (approximate) solution to the reduced
problem, it recovers an (approximate) solution to the original problem, \ie,
$(x^\star, y^\star, z^\star)$ (approximately) satisfying \eqref{eq:KKT}.

\section{Experiments}
\label{sec:experiments}
In this section we first compare the presolve performance of PSLP against the
commercial state-of-the-art presolver implemented in Gurobi
\cite{Achterberg2019}. There are two main reasons for this comparison. First, we
wish to gain insight into the relative importance of fast and medium explorers
versus slow explorers, as PSLP implements only the former while Gurobi's
presolver implements all three categories. Second, since a central goal of this
paper is to assess whether presolving improves the performance of
GPU-accelerated first-order LP solvers, it is essential that the presolver
itself is reasonably effective. An inefficient presolver would prevent a fair
and honest assessment of the impact of presolving and could distort the results.
Comparing PSLP to a commercial state-of-the-art presolver therefore also serves
to validate that it is adequate for a fair assessment of the impact of presolving.

After comparing PSLP to Gurobi's presolver, we evaluate how presolving using PSLP
or Gurobi affects the performance of cuPDLPx \cite{Lu2025-cuPDLPx}, a
GPU-accelerated first-order LP solver based on PDHG.

\paragraph{Benchmark datasets.} We conduct numerical experiments on two LP
benchmark datasets. First, we consider Mittelmann's LP collection
\cite{mittelmannLP}, which consists of 48 publicly available instances with more
than $10^5$ nonzero entries. Second, we consider 383 instances derived from
root-node LP relaxations of MIPLIB 2017 \cite{Gleixner2021}, using the filtering
criteria described in \cite{Lu2025}. We categorize problems into three groups
based on their size: \emph{small} ($10^5$-$10^6$ nonzero entries), \emph{medium}
($10^6$-$10^7$ nonzero entries), and \emph{large} (more than $10^7$ nonzero
entries).

\paragraph{Hardware.} Problems are presolved on the CPU and solved on 
the GPU. For the presolve step, we use a Macbook M4 pro (14 cores)
with 48 GB of unified memory, and then execute cuPDLPx on a 
NVIDIA H100-SXM-80GB GPU with CUDA 12.4. 

\subsection{Presolve performance}

\paragraph{Setup.} We compare PSLP with Gurobi's presolver (version 13.0, run
with default settings) using the following metrics: 
\begin{itemize}
    \item Arithmetic and geometric mean of the presolve time, given in seconds
    (denoted by ``AM-T'' and ``GM-T'', respectively). For each problem instance,
    we average the presolve time over 10 runs.
    \item The average reduction in the number of nonzero entries, given as the
    ratio of the number of nonzero entries in the reduced problem to that in the
    original problem (denoted by ``ratio (nnz)''). 
\end{itemize}

\paragraph{Results.}
Table \ref{tab:presolver_comparison} summarizes the results. On Mittelmann's LP
collection, the average per-instance ratio of Gurobi's presolve time to PSLP's
presolve time is 11.8 (not shown in the table). Despite this substantial speedup
of an order of magnitude, PSLP on average achieves 90\% of the reduction in the
number of nonzero entries achieved by Gurobi's presolver. On the MIPLIB
relaxations, the corresponding average time ratio is 6.58 (also not shown in the
table), while PSLP on average attains 94\% of the reduction in the number of nonzero
entries achieved by Gurobi's presolver.

\begin{figure}
    \centering
    \includegraphics[width=0.85\textwidth]{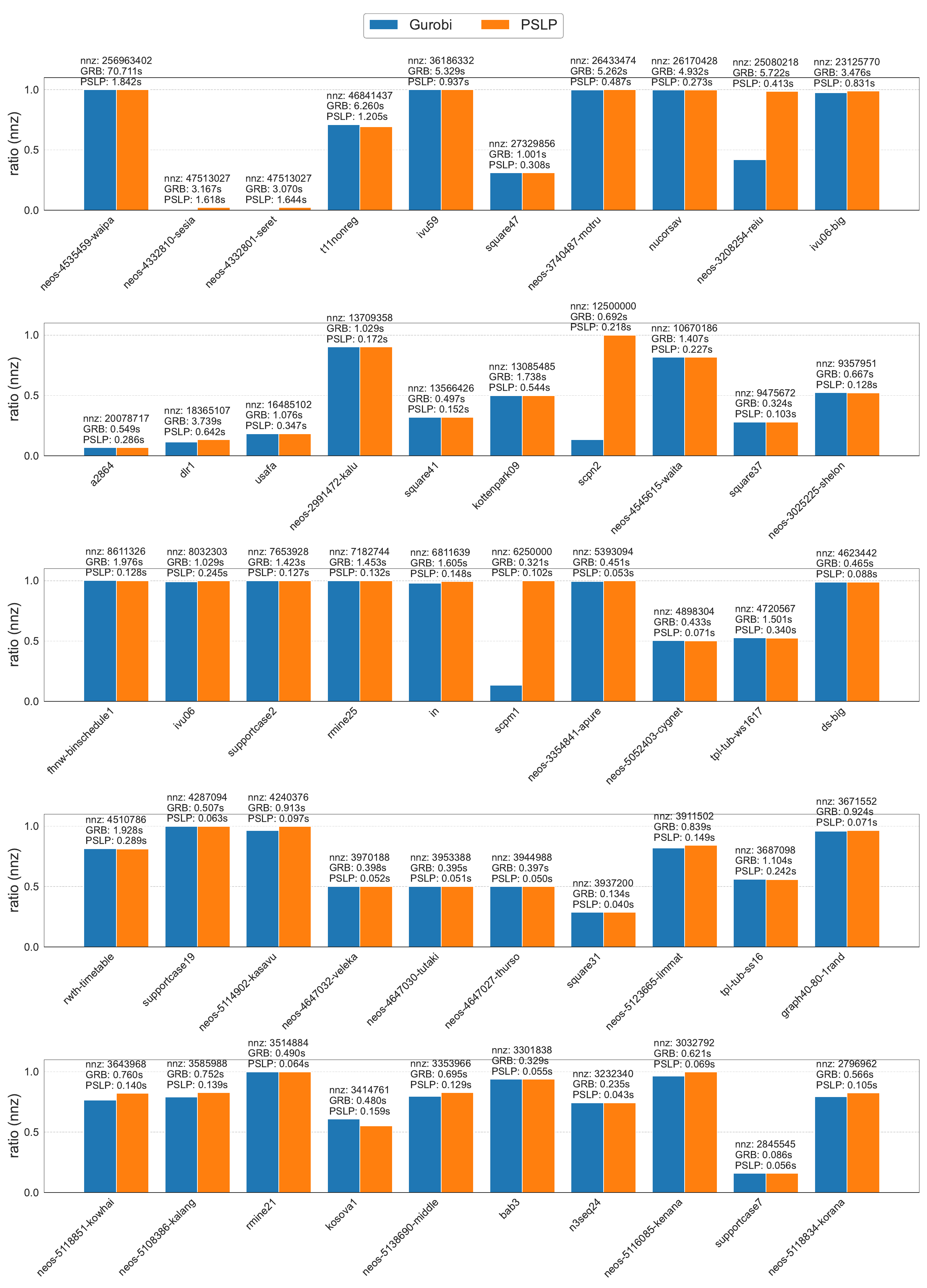}
    \caption{Comparison of PSLP and Gurobi's presolver on the 50 largest MIPLIB relaxations.}
    \label{fig:pslp_vs_gurobi_miplib}
\end{figure}

Figures \ref{fig:pslp_vs_gurobi_mittelmann} and \ref{fig:pslp_vs_gurobi_miplib}
compare the reductions in the number of nonzero entries achieved by PSLP and
Gurobi's presolver across the two datasets, together with the corresponding
presolve times. (For the MIPLIB relaxations, results are reported only for the
50 largest instances due to space constraints.) Overall, PSLP achieves
reductions comparable to those of Gurobi's presolver on all but five instances
(\verb|L1_sixm1000obs|, \verb|scpm1|, and \verb|L1_sixm250obs| from Mittelmann's
LP collection, and \verb|neos-3208254-reui| and \verb|scpm2| from the MIPLIB
relaxations), while often requiring significantly less presolve time.

\begin{table}
\centering
\small
\begin{tabular}{llccc}
\toprule
\textbf{Dataset} & \textbf{Presolver} & \textbf{AM} & \textbf{GM} & \textbf{ratio (nnz)} \\
\midrule
\multirow{2}{*}{Mittelmann}
 & Gurobi & 3.59 & 0.70  & 0.71 \\
 & PSLP   & 0.33 & 0.084 & 0.79 \\
\midrule
\multirow{2}{*}{MIPLIB relaxations}
 & Gurobi & 0.47 & 0.077 & 0.76 \\
 & PSLP   & 0.057 & 0.013  & 0.81 \\
\bottomrule
\end{tabular}
\caption{Comparison of PSLP and Gurobi's presolver. The arithmetic mean (AM) and 
geometric mean (GM) of the presolve time are given in seconds.}
\label{tab:presolver_comparison}
\end{table}

\begin{figure}
    \centering
    \includegraphics[width=0.85\textwidth]{figures/reduction_histogram_mittelmann.pdf}
    \caption{Comparison of PSLP and Gurobi's presolver on Mittelmann's LP
    collection. Each bin shows the ratio of the nonzero entries in the reduced
    problem to that in the original problem. Above each bin, we show the number
    of nonzero entries in the original problem, and the presolve times of PSLP
    and Gurobi. The problems are sorted by the number of nonzero
    entries in the original problem.}
    \label{fig:pslp_vs_gurobi_mittelmann}
\end{figure}

\subsection{End-to-end solver performance: experimental setup}

\paragraph{Solver.} We run cuPDLPx (version 0.2.5) both with and without PSLP
and Gurobi's presolver enabled. The termination criterion adopted by cuPDLPx
depends on a relative tolerance parameter $\epsilon_{\text{rel}} > 0$. We use
$\epsilon_{\text{rel}} = 10^{-8}$ in all experiments, which corresponds to the
default value for higher accuracy \cite{Lu2025-cuPDLPx}.

\paragraph{Time limit.} For the MIPLIB relaxations, we set a time limit of 3600
seconds for small and medium-sized problems, and 18000 seconds for large
problems. For Mittelmann's LP collection, we set a time limit of 15000 seconds.
(These time limits are the same as those used in \cite{Lu2025-cuPDLPx}.) If an
instance remains unsolved when the time limit is reached, we count it as a
failure and set its solve time to the time limit for the purpose of computing
average solve times.

\paragraph{Performance metrics.} We consider the following common metrics: 
\begin{itemize}
\item \textbf{Count}: The number of problems solved within the time limit.
\item \textbf{SGM}: We use the shifted geometric mean of the total solve time,
given in seconds. The shifted geometric mean is defined as $(\prod_{i=1}^K(t_i +
\Delta))^{1/K} - \Delta$, where $t_i$ is the total solve time for the $i$th
instance (in seconds), $K$ is the total number of instances, and $\Delta > 0$ is
the shift. We shift by $\Delta = 1$ and $\Delta = 10$ seconds and denote these
metrics as SGM1 and SGM10, respectively. Without presolving, the total solve
time is simply the time taken by cuPDLPx to solve the original problem. With
presolving, the total solve time is the sum of the presolve time and the time
taken by cuPDLPx to solve the reduced problem.
\item \textbf{Win rate}: The percentage of problem instances for which the total
solve time with presolve enabled is less than the total solve time with presolve
disabled. 
\end{itemize}

\subsection{End-to-end solver performance: results}

\begin{figure}
    \centering
    \begin{subfigure}[b]{0.95\textwidth}
        \centering
        \includegraphics[width=\textwidth]{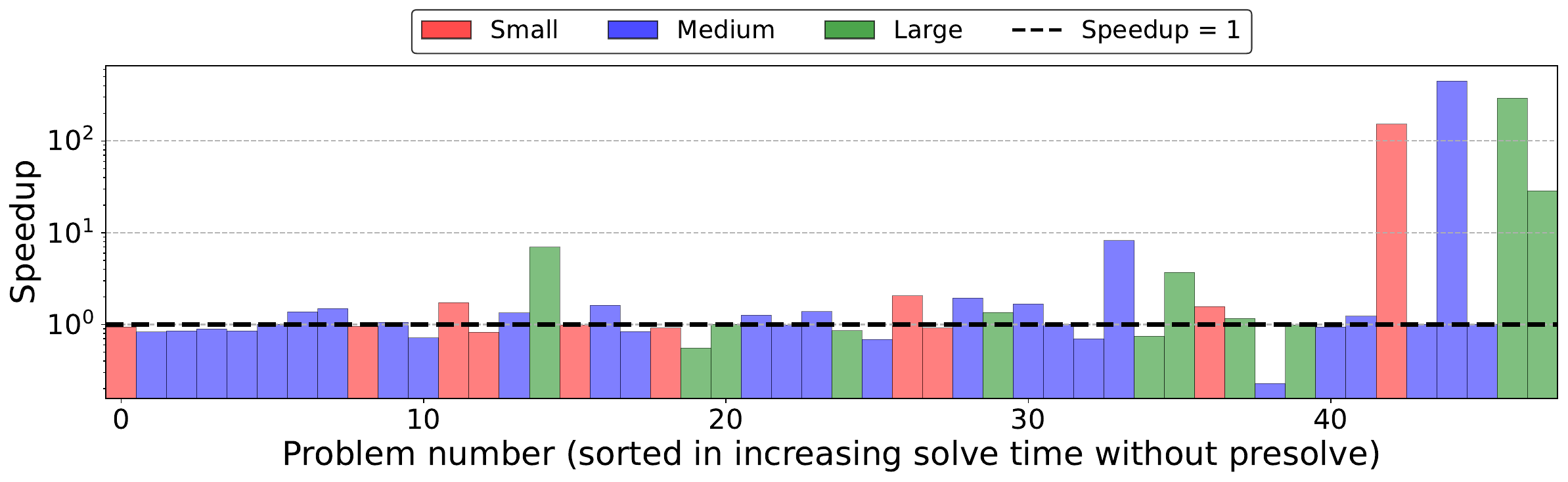}
        \caption{Mittelmann's LP collection.}
    \end{subfigure}
    \vspace{1em}
    \begin{subfigure}[b]{0.95\textwidth}
        \centering
        \includegraphics[width=\textwidth]{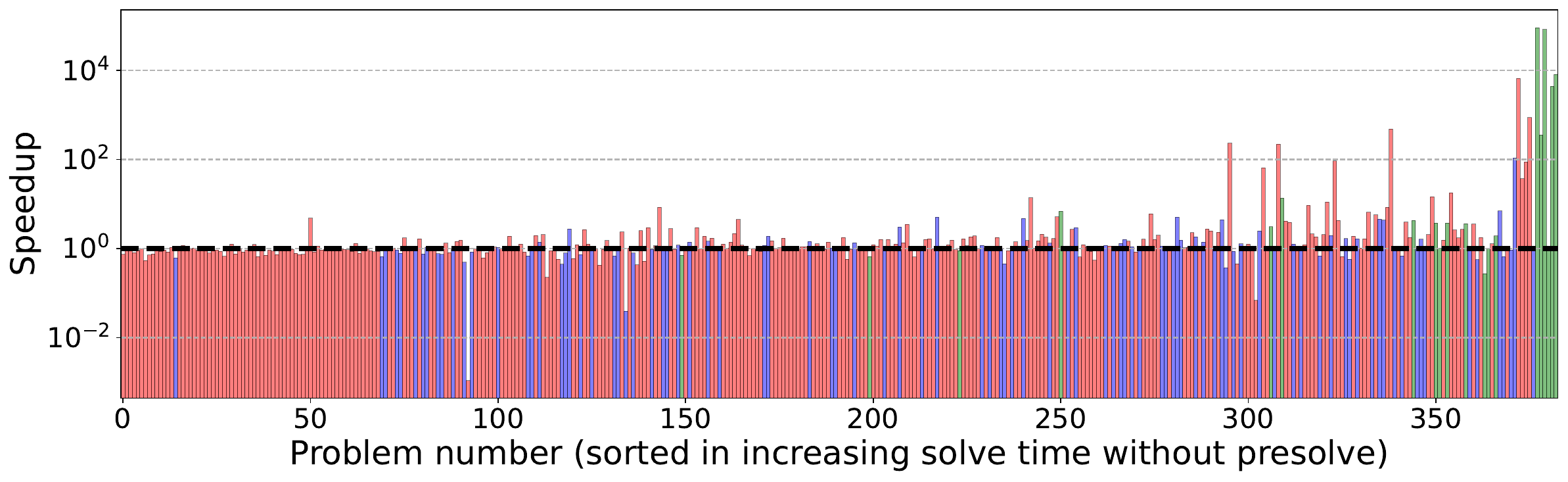}
        \caption{MIPLIB relaxations.}
    \end{subfigure}
    \caption{The speedup (\ie, the ratio of the total solve time without
    presolve to the total solve time with presolve) of enabling PSLP for cuPDLPx
    for each problem instance.}
    \label{fig:speedup_histogram}
\end{figure}

\paragraph{Presolve impact on total solve time.}
Table \ref{tab:results} summarizes the results. We see that enabling PSLP
improves the shifted geometric mean in all three categories for both datasets
compared to the case with presolve disabled. The improvement is more significant
for the large problems, with almost a three-fold improvement for Mittelmann's LP
collection and more than a nine-fold improvement for the MIPLIB relaxations. We
also see that the solve times with either PSLP or Gurobi's presolver enabled are
similar. Interestingly, given the significant improvement in solve time with
presolve enabled, the win rate (\ie, the percentage of instances for which
presolve enabled is faster than presolve disabled) is lower than one might
expect. One possible explanation is that for problems where cuPDLPx already
converges fast with a modest number of iterations, presolving becomes harmful.
On the other hand, there exist problems where cuPDLPx is very slow without
presolve but becomes fast with presolve, which gives a significant improvement
in solve time. We illustrate this in Figure \ref{fig:speedup_histogram}, which
for each problem instance shows the speedup of enabling PSLP for cuPDLPx. In the
right part of the figure, which corresponds to problems with the longest solve
time without presolve, we see that the speedup from presolve can be orders of
magnitude. 

\begin{table}
\centering
{\footnotesize
\setlength{\tabcolsep}{5.3pt}
\begin{tabular}{lcrrc|crrc}
\toprule
& \multicolumn{4}{c|}{\textbf{Mittelmann}} & \multicolumn{4}{c}{\textbf{MIPLIB}} \\
& \textbf{Count} & \textbf{SGM1} & \textbf{SGM10} & \textbf{WR (\%)} & \textbf{Count} & \textbf{SGM1} & \textbf{SGM10} & \textbf{WR (\%)} \\ \midrule
\textbf{Small problems} & 10 & ------ & ------ & ------ & 269 & ------ & ------ & ------ \\
\hline
cuPDLPx         & 9  & 16.23 & 33.20 &  ------ & 265  & 2.27 & 6.08 & ------ \\
cuPDLPx + PSLP  & 10  & 8.31 & 14.23 & 40 & 269 & 1.48 & 3.54 & 60 \\
cuPDLPx + Gurobi   & 10 & 11.83 & 19.75 &  50 &  269 & 1.40 & 3.01 & 44  \\ \midrule
\textbf{Medium problems} & 27 & ------ & ------ & ------ & 94 & ------ & ------ & ------ \\
\hline
cuPDLPx         & 24  & 28.85 & 60.76 & ------  & 92  & 5.57 & 12.28 & ------ \\
cuPDLPx + PSLP  & 25  & 21.15 & 45.17 & 52 &  93 & 4.62 & 10.54 & 53 \\
cuPDLPx + Gurobi   &  25  & 19.70 & 39.26 & 41  & 93 & 4.87 & 10.01 & 34 \\ \midrule
\textbf{Large problems} & 11 & ------ & ------ & ------ & 20 & ------ & ------ & ------ \\
\hline
cuPDLPx         & 9  & 136.51 & 169.33 & ------ & 14  & 286.03 & 408.26 & ------  \\
cuPDLPx + PSLP  & 11 & 50.44 & 63.05 & 64 &  19 & 22.12 & 44.12 & 75 \\
cuPDLPx + Gurobi   & 11  & 56.10 & 68.05 & 36 & 19 & 22.97 & 42.18 & 70  \\
\midrule
\textbf{All problems} & 48 & ------ & ------ & ------ & 383 & ------ & ------ & ------ \\
\hline
cuPDLPx         & 42  & 36.78 & 69.01 & ------ &  371 & 3.90 & 10.77 & ------  \\
cuPDLPx + PSLP  & 46 & 21.43 & 39.57 & 52 &   381 & 2.41 & 6.13 & 59 \\
cuPDLPx + Gurobi   & 46  & 22.64 & 39.35 & 42 & 381 & 2.37 & 5.54 & 43  \\
\bottomrule
\end{tabular}
}
\caption{Results on Mittelmann's LP collection and MIPLIB relaxations. The
abbreviation WR stands for win rate, which is the percentage of instances for
which presolve enabled is faster than presolve disabled, and SGM1 and SGM10 are given in
seconds.}
\label{tab:results}
\end{table}

\paragraph{Time spent in presolve.}
We analyze the ratio between the time spent in presolve and the time spent on
solving the \emph{reduced} problem. Table \ref{tab:presolve_time} summarizes the
average ratio across different problem categories. We see that the presolve time
for PSLP is generally a small fraction of the solve time that tends to be more
significant for larger problems. Gurobi's presolve time is even more
significant and sometimes dominates the time spent on solving the reduced
problem.

\begin{table}[h]
\centering
{\small
\begin{tabular}{lcccc|cccc}
\toprule
& \multicolumn{4}{c|}{\textbf{Mittelmann}} & \multicolumn{4}{c}{\textbf{MIPLIB}}
\\
& \textbf{Small} & \textbf{Medium} & \textbf{Large} & \textbf{All} &
\textbf{Small} & \textbf{Medium} & \textbf{Large} & \textbf{All} \\ \midrule
PSLP     (\%)   & 1.3 & 12.6  & 19.3 & 11.8 & 7.4   & 20.6  & 39.1  & 12.1  \\
Gurobi   (\%)   & 6.2 & 100.9 & 82.3 & 76.9 & 29.1  & 105.9 & 198.3 & 56.0  \\
\bottomrule
\end{tabular}
} \caption{Average presolve time as a percentage of the time it takes for
cuPDLPx to solve the reduced problem.}
\label{tab:presolve_time}
\end{table}

\section{Discussion}
\label{sec:discussion}
We have revisited LP presolving, a classical topic in the optimization
literature, from the perspective of GPU-accelerated first-order solvers. While a
myriad of presolve techniques have been proposed over the decades, a systematic
investigation of their cost-effectiveness in the context of GPU-accelerated
first-order methods has been lacking. To address this, we organized these
techniques into a three-class taxonomy based on computational cost and
implemented the two cheaper classes in PSLP, a new open-source presolver
designed to be lightweight. By comparing PSLP to the commercial state-of-the-art
presolver in Gurobi, which implements all three classes of presolve techniques,
our experiments on standard benchmarks reveal several findings:
\begin{itemize}
    \item A lightweight presolver can achieve reductions comparable to those of
    a sophisticated state-of-the-art commercial presolver, at a fraction of the
    computational cost (see table \ref{tab:presolver_comparison} and figures
    \ref{fig:pslp_vs_gurobi_mittelmann} and \ref{fig:pslp_vs_gurobi_miplib}).
    \item The time spent in presolve can be significant relative to the time
    spent on solving the reduced problem (see table \ref{tab:presolve_time}).
    \item Lightweight presolve can substantially improve the end-to-end
    performance of GPU-accelerated first-order solvers (see table
    \ref{tab:results}). While many problems are largely unaffected by presolve,
    the speedup can be orders of magnitude for certain larger problems
    (see figure \ref{fig:speedup_histogram}).
\end{itemize}
Presolve has traditionally been regarded as a preprocessing step whose cost is
negligible relative to the solve phase. Our first two findings challenge this
assumption in the context of GPU-accelerated first-order methods, and together
highlight the growing importance of presolve cost. As GPU hardware improves, the
solve phase of first-order methods will continue to accelerate, but the
sequential nature of presolving means it will not benefit from these advances to
the same extent. Over time, presolve therefore risks becoming an increasingly
dominant fraction of the total solve time. Our results show that by restricting
to cheaper presolve techniques, much of the reduction benefit can be retained
without presolve becoming an excessive bottleneck. Nonetheless, even with a
lightweight presolver such as PSLP, the presolve time is already non-negligible
(see table \ref{tab:presolve_time}). As first-order solvers become faster, the
overall speedup in total solve time may therefore be limited unless the presolve
phase itself is accelerated. (For example, if cuPDLPx became twice as fast due
to hardware improvements, a rough estimate based on table
\ref{tab:presolve_time} suggests that the speedup in total solve time for large
MIPLIB problems would be about $1.20\times$ if Gurobi's presolver were used,
rather than the expected $2\times$.) While skipping presolve altogether is
possible, our results demonstrate the benefit of presolving with PSLP,
suggesting that it should generally not be skipped in solvers similar to
cuPDLPx. More broadly, future research should consider accelerating not only the
solve phase but also the presolve phase, which has been largely neglected in the
literature. Given the substantial engineering effort required to develop an
efficient presolver, such research questions have until now been impractical to
explore. By open-sourcing PSLP, we hope to enable and encourage research in this
direction. 

Our third finding highlights the continued importance of presolving for linear
programming, even in the current era of GPU-accelerated first-order solvers. It
also indicates that presolving has a stronger impact on first-order solvers (at
least cuPDLPx) than is suggested by the seminal work of
\cite{Applegate2021}, where PaPILO \cite{Gleixner2023} was used for presolving
and only very marginal reductions in iteration counts were reported
when presolving was applied prior to a first-order method. Furthermore, the
presolve time (which is significant relative to the solve time for the reduced
problem) was not included in the performance metrics reported in
\cite{Applegate2021}.

An efficient LP presolver forms the foundation of presolvers for more general
problem classes. In ongoing work, we are extending PSLP to support quadratic
programs and convex conic programs involving second-order cones and exponential
cones.

\section*{Acknowledgments}
We are very grateful to Zedong Peng, main developer of cuPDLPx, for collecting
the data for the cuPDLPx experiments in \S4, and for providing exceptionally
detailed and valuable feedback on PSLP. Without his input, PSLP would not be in
the shape it is. We also thank Max Schaller, Kasper Johansson, Chris Maes,
and Rajesh Gandham for helpful feedback and discussions.

\newpage
\bibliographystyle{alpha}
\bibliography{references.bib}
\end{document}